\documentclass[12pt,reqno]{amsart}

\usepackage{amsmath,amsthm,amsfonts,amssymb}
\usepackage{cite}
\usepackage[dvipsnames]{xcolor}
\usepackage[pdfpagelabels,colorlinks=true,linkcolor=blue,citecolor=OliveGreen,urlcolor=blue]{hyperref}

\newtheorem{theorem}{Theorem}[section]
\newtheorem{lemma}[theorem]{Lemma}
\newtheorem{proposition}[theorem]{Proposition}
\theoremstyle{definition}
\newtheorem{definition}[theorem]{Definition}
\theoremstyle{remark}

\numberwithin{equation}{section}

\DeclareMathOperator*{\esssup}{ess\,sup}
\DeclareMathOperator*{\essinf}{ess\,inf}
\newcommand{\R}{\mathbb{R}}
\begin{document}

\title[Blow-up time estimates]{Blow-up and blow-up-time estimates for a singular pseudo-parabolic equation with a space--time variable exponent}

\author{Nguyen Thanh Tung}
\address{Nguyen Thanh Tung, Ho Chi Minh City University of Education, Vietnam}
\email{tungnthanh@hcmue.edu.vn}

\author{Le Xuan Truong}
\address{Le Xuan Truong, University of Economics Ho Chi Minh City (UEH), Vietnam}
\email{lxuantruong@ueh.edu.vn}

\author{Tan Duc Do}
\address{Tan Duc Do, University of Economics Ho Chi Minh City (UEH), Vietnam}
\email{tandd.am@ueh.edu.vn}

\author{Nguyen Ngoc Trong}
\address{Nguyen Ngoc Trong, Group of Analysis and Applied Mathematics, Department of Mathematics, Ho Chi Minh City University of Education, Vietnam}
\email{trongnn@hcmue.edu.vn}
\thanks{Nguyen Ngoc Trong is the corresponding author.}


\begin{abstract}
We consider a singular pseudo-parabolic initial-boundary value problem. Its time derivative carries the factor $|x|^{-2}$, while the source exponent $p(x,t)$ depends on both space and time. This time dependence produces a logarithmic term in the energy balance, so the usual autonomous energy is no longer monotone. We introduce a modified energy that compensates for this term and obtain finite-time blow-up criteria in the negative- and nonnegative-energy regimes, together with explicit upper bounds for the blow-up time. We also derive a lower bound for every finite blow-up time by a direct Sobolev estimate. The arguments require only the subcritical condition $2<p^-\le p(x,t)\le p^+<2^*$, under the stated regularity, monotonicity, and boundedness assumptions on $p$ and $k$.
\end{abstract}

\maketitle

\section{Introduction}

Variable-exponent problems arise naturally in models whose constitutive response changes with position or time. Classical examples include nonlinear elasticity and electrorheological fluids \cite{Hal,AM,DR,Ruz}, while variable-growth diffusion has also been used in image processing \cite{ANA}. The functional-analytic framework for Lebesgue and Sobolev spaces with variable exponents is by now well developed; we refer to \cite{DHHR} and the references therein.

Finite-time blow-up is another basic theme in nonlinear evolution equations. Depending on the structure of the equation, blow-up can be detected through eigenfunction arguments, energy methods, potential wells, concavity methods, or comparison principles. General background can be found in \cite{Hu,GV,Lev}. For equations closer to the present setting, lower and upper estimates of the blow-up time have been studied for semilinear parabolic equations, pseudo-parabolic equations, and higher-order diffusion models; see, for example, \cite{BGH,SLW,TTDa,TTDb}.

Two lines of work are especially close to the present problem. On the variable-exponent side, time-dependent exponents in nonlinear parabolic equations were studied in \cite{GG}, pseudo-parabolic equations with variable exponents were considered in \cite{DSP}, and variable-exponent reaction--diffusion equations with a special medium void and damping effects were studied in \cite{Do}. On the singular pseudo-parabolic side, the constant-exponent problem
\[
\frac{u_t}{|x|^s}-\Delta u_t-\Delta u=|u|^{p-2}u,
\qquad 0\le s\le2,
\]
was studied systematically by Lian, Wang and Xu \cite{LWX}, who obtained global-existence and blow-up results together with upper bounds for the blow-up time. Subsequent work gave further blow-up-time and blow-up-rate estimates \cite{SHO} and new blow-up/lower-bound results \cite{YPC}. Earlier singular or degenerate filtration models associated with a special medium void include \cite{Tan,Han,Han2}.

Accordingly, the novelty here is not the study of blow-up-time bounds for singular pseudo-parabolic equations by itself. Rather, we combine the singular pseudo-parabolic structure with a genuinely space--time dependent source exponent $p(x,t)$ and a nonconstant coefficient $k(t)$. This changes the energy mechanism in an essential way: differentiating the nonlinear potential produces an additional $p_t\log|u|$ term that is absent from the constant-exponent models \cite{LWX,SHO,YPC}. The modified energy introduced below is designed precisely to compensate for this term.

The time dependence of $p(x,t)$ creates the principal technical difficulty. Differentiating the nonlinear potential
\[
\int_\Omega \frac{|u(x,t)|^{p(x,t)}}{p(x,t)}\,dx
\]
produces, besides the usual duality term with $u_t$, a logarithmic contribution containing
\[
\frac{p_t(x,t)}{p(x,t)^2}
\bigl[p(x,t)\log |u(x,t)|-1\bigr]|u(x,t)|^{p(x,t)}.
\]
This term has no fixed sign, so the standard autonomous energy is not directly monotone. A second issue is that the singular kinetic term naturally leads to the Hilbert quantity
\[
\frac12\left(\left\|\frac{u}{|x|}\right\|_2^2+\|\nabla u\|_2^2\right),
\]
and both components must be retained in the concavity argument.

Our approach is to add a time-dependent correction to the usual energy. Under $p_t\ge0$ and a nondecreasing bounded coefficient $k$, this modified energy yields a monotone quantity that controls the logarithmic term. We then obtain three results. For negative modified initial energy, a differential inequality gives finite-time blow-up and an explicit upper bound. For nonnegative modified energy below a computable threshold, a corrected Levine-type concavity argument again yields finite-time blow-up and an upper bound. Finally, if a solution blows up, a direct Sobolev estimate gives an integral lower bound for the blow-up time. In particular, the lower-bound argument does not require an additional Gagliardo--Nirenberg restriction beyond the Sobolev-subcritical assumption.

We emphasize that after placing $u_t$ in the natural $W^{1,2}_0(\Omega)$ class, the analysis needs only
\[
2<p^-\le p(x,t)\le p^+<2^*:=\frac{2n}{n-2}.
\]
This is the exponent range used throughout the revised arguments below.

Specifically, let $n \in \{3,4,5,\ldots\}$.
Let $\Omega \subset \R^n$ be open bounded with Lipschitz boundary.
Denote $Q := \Omega \times (0,\infty)$.
Let $p: \overline{Q} \longrightarrow (0,\infty)$ satisfy 
\begin{equation} \label{eq:p-assumption}
	\begin{cases}
		p \in C^1(\overline{Q}), \quad p_t \ge 0,
		\\[2mm]
		\displaystyle
		2 < p^- \le p(\cdot) \le p^+ < 2^*:=\frac{2n}{n-2}.
	\end{cases}
\end{equation}
Here
\[
p^- := \essinf_{(x,t) \in Q} p(x,t)
\quad \mbox{and} \quad
p^+ := \esssup_{(x,t) \in Q} p(x,t).
\]
Next let $k$ be a function with the following properties:
\begin{equation} \label{eq:k-assumption}
	\begin{cases}
		k \in C^1[0,\infty), 
		\quad k(0) > 0
		\quad \mbox{and} \quad
		k'(t) \ge 0  \,\, \mbox{ for all } t \in [0,\infty),
		\\[2mm]
		k_\infty := \displaystyle\lim_{t\to\infty} k(t) < \infty.
	\end{cases}
\end{equation}

In this paper, we consider the following pseudo-parabolic equation with singular potential:
\[
(P) \quad	\left\{\begin{array}{ll}
	\displaystyle\frac{u_t}{|x|^2} - \Delta u_t- \Delta u = k(t) \, |u|^{p(x,t)-2}u, & (x,t) \in \Omega_T := \Omega \times (0,T),
	\\
	u(x,t) = 0, & (x,t) \in \partial \Omega \times (0,T), \smallskip\\
	u(x,0) = u_0(x), & x \in \Omega,
\end{array}\right.
\]
where $T > 0$ and $0 \ne u_0 \in W^{1,2}_0(\Omega)$.
When $p$ is constant and $k\equiv1$, this is the critical case $s=2$ of the singular pseudo-parabolic family studied in \cite{LWX}; see also \cite{SHO,YPC} for subsequent blow-up-time results. If the singular factor $|x|^{-2}$ is removed from the time derivative, $p(x,t)$ is replaced by a spatial exponent, and $k\equiv1$, one obtains the type of generalized heat equation studied in \cite{AMK}. Related blow-up models can be found in \cite{BGH,DSP,Do,SLW,TTDa,TTDb}.

We now state the solution concept and the natural energy quantity used throughout.

\begin{definition}
	Let $p$ be given by \eqref{eq:p-assumption} and $0 \ne u_0 \in W^{1,2}_0(\Omega)$.
	A function $u$ is called a weak solution to $(P)$ on $[0,T]$ if
	\[
		u \in C([0,T];W^{1,2}_0(\Omega)),
		\qquad
		u_t \in L^2(0,T;W^{1,2}_0(\Omega)),
	\]
	$u(0)=u_0$, and, for every $\varphi\in W^{1,2}_0(\Omega)$ and for a.e.\ $t\in(0,T)$,
	\begin{equation} \label{eq:weak-form}
		\int_\Omega \frac{u_t\varphi}{|x|^2}\,dx
		+ (\nabla u_t,\nabla\varphi)
		+ (\nabla u,\nabla\varphi)
		= k(t)\int_\Omega |u|^{p(x,t)-2}u\,\varphi\,dx.
	\end{equation}
\end{definition}
By Hardy's inequality, the regularity $u_t\in L^2(0,T;W^{1,2}_0(\Omega))$ automatically implies
\[
\int_0^T\int_\Omega \frac{|u_t(x,t)|^2}{|x|^2}\,dx\,dt<\infty,
\]
so the singular term in \eqref{eq:weak-form} is well defined.

For a weak solution $u$, we use throughout the natural energy quantity
\begin{equation} \label{eq:L-def}
	L(t):=\frac12\left(
	\left\|\frac{u(t)}{|x|}\right\|_{L^2(\Omega)}^2
	+\|\nabla u(t)\|_{L^2(\Omega)}^2
	\right).
\end{equation}

\begin{definition}
	Let $u$ be a weak solution to $(P)$ on its maximal interval of existence $[0,T^*)$, where $T^*\in(0,\infty]$.
	If $T^*=\infty$, then $u$ is called a global solution. If $T^*<\infty$, we say that $u$ blows up at $T^*$ when
	\begin{equation} \label{eq:blowup}
		\limsup_{t\to(T^*)^-}L(t)=\infty.
	\end{equation}
\end{definition}

For $u\in W^{1,2}_0(\Omega)$ and $t\ge0$, define the energy and Nehari functionals by
\begin{align}
J(u,t)
&:=\frac12\int_\Omega |\nabla u|^2\,dx
-k(t)\int_\Omega \frac{|u|^{p(x,t)}}{p(x,t)}\,dx,
\label{eq:J-def}\\
I(u,t)
&:=\int_\Omega |\nabla u|^2\,dx
-k(t)\int_\Omega |u|^{p(x,t)}\,dx.
\label{eq:I-def}
\end{align}
Because the exponent depends on time, it is convenient to introduce the modified energy
\begin{equation}\label{eq:E-def}
E(u,t):=J(u,t)+k_\infty\int_\Omega \frac{1}{p(x,t)}\,dx.
\end{equation}
The correction term in \eqref{eq:E-def} is designed to control the logarithmic contribution generated by differentiating $|u|^{p(x,t)}$ with respect to time.

Our first main result concerns an upper bound on the blow-up time for a weak solution to $(P)$ when the modified initial energy is negative.

\begin{theorem} \label{thm:negative}
	Let $n \in \{3,4,5,\ldots\}$ and $\Omega \subset \R^n$ be open bounded with Lipschitz boundary. 
	Let $p$, $k$ satisfy \eqref{eq:p-assumption} and \eqref{eq:k-assumption} respectively.
	Suppose further that 
	\[
	E(u_0,0) < 0.
	\]
	Let $u$ be the maximal weak solution to $(P)$.
	Then $u$ blows up at a finite time $T^*$ satisfying
	\[
	T^* 
	\le \frac{\left\| \displaystyle\frac{u_0}{|x|} \right\|_{L^2(\Omega)}^2+\|\nabla u_0\|_{L^2(\Omega)}^2}
	{p^- \, (p^- - 2) \, [-E(u_0,0)]}.
	\]
\end{theorem}

A second upper bound is available when the modified initial energy is nonnegative but lies below a quantitative threshold.

\begin{theorem} \label{thm:positive}
	Let $n \in \{3,4,5,\ldots\}$ and $\Omega \subset \R^n$ be open bounded with Lipschitz boundary.
	Let $p$, $k$ satisfy \eqref{eq:p-assumption} and \eqref{eq:k-assumption} respectively, and set
	\[
	C_1:=\frac{p^-\,(1+H_n)}{p^- -2},
	\qquad H_n:=\frac{4}{(n-2)^2}.
	\]
	Suppose that the initial modified energy is nonnegative and satisfies
	\[
	0\le C_1E(u_0,0)<L(0),
	\qquad
	L(0)=\frac12\left(
	\left\|\frac{u_0}{|x|}\right\|_{L^2(\Omega)}^2
	+\|\nabla u_0\|_{L^2(\Omega)}^2
	\right).
	\]
	Let $u$ be the maximal weak solution to $(P)$ and define
	\[
	M_0:=L(0)-C_1E(u_0,0)>0.
	\]
	Then $u$ blows up at a finite time $T^*$ satisfying
	\[
	T^*\le
	\frac{8C_1(p^- -1)L(0)}
	{p^-(p^- -2)^2M_0}.
	\]
\end{theorem}

Lastly we present a lower bound on the blow-up time.

\begin{theorem} \label{thm:lower}
	Let $n \in \{3,4,5,\ldots\}$ and $\Omega \subset \R^n$ be open bounded with Lipschitz boundary.
	Let $p$, $k$ satisfy \eqref{eq:p-assumption} and \eqref{eq:k-assumption}, respectively, and put
	\[
	\gamma^-:=\frac{p^-}{2},
	\qquad
	\gamma^+:=\frac{p^+}{2}.
	\]
	Let $u$ be a maximal weak solution to $(P)$ which blows up at a finite time $T^*$.
	Then there exists a constant
	$C^*=C^*(\Omega,n,p^\pm,k_\infty)>0$ such that, for every $t_0\in[0,T^*)$,
	\[
	T^* \ge t_0+\frac{1}{C^*}
	\int_{L(t_0)}^\infty
	\frac{ds}{s^{\gamma^-}+s^{\gamma^+}}.
	\]
	In particular, if $L(t_0)\ge1$, then
	\[
	T^*\ge t_0+
	\frac{L(t_0)^{1-\gamma^+}}
	{2C^*(\gamma^+-1)}.
	\]
\end{theorem}

The paper is organized as follows.
Section \ref{sec:prelim} collects the Hardy and Sobolev inequalities together with the energy identities used throughout the paper.
The upper and lower bounds on the blow-up time are established in Sections \ref{sec:upper} and \ref{sec:lower}, respectively.

\section{Preliminaries} \label{sec:prelim}

We collect only the estimates needed in the blow-up arguments. The assumptions in \eqref{eq:p-assumption} ensure throughout that the nonlinear integrals below are finite by the standard Sobolev embedding $W^{1,2}_0(\Omega)\hookrightarrow L^q(\Omega)$ for $2\le q<2^*$.

\subsection{Fundamental inequalities}

We begin with the Hardy inequality associated with the singular factor $|x|^{-2}$.

\begin{lemma}
	Let $n \ge 3$ and $u \in W^{1,2}_0(\Omega)$.
	Then $\displaystyle\frac{u}{|x|} \in L^2(\Omega)$ and
	\[
	\int_{\Omega} \frac{|u|^2}{|x|^2} \, dx 
	\le \frac{4}{(n-2)^2} \, \int_{\Omega} |\nabla u|^2 \, dx
	=: H_n \, \int_{\Omega} |\nabla u|^2 \, dx.
	\]
\end{lemma}

\begin{proof}
	This follows at once from \cite[Theorem 4.2.2]{BEL}.
\end{proof}

We shall also use the standard subcritical Sobolev embedding.

\begin{lemma}\label{lem:sobolev}
	Let $n \ge 3$ and $2<q<2^*$.
	Then there exists a constant $S_q=S_q(\Omega,n,q)>0$ such that
	\[
	\|u\|_{L^q(\Omega)} \le S_q\|\nabla u\|_{L^2(\Omega)}
	\]
	for every $u\in W^{1,2}_0(\Omega)$.
\end{lemma}

The following standard local well-posedness observation records the continuation property needed in the blow-up arguments.

\begin{proposition}\label{prop:local}
For every $u_0\in W^{1,2}_0(\Omega)$ there exists a unique maximal weak solution $u$ of $(P)$ on an interval $[0,T^*)$, $T^*\in(0,\infty]$. Moreover, if $T^*<\infty$, then
\[
\limsup_{t\to(T^*)^-}L(t)=\infty.
\]
\end{proposition}

\begin{proof}
Set $H:=W^{1,2}_0(\Omega)$ with norm $\|v\|_H:=\|\nabla v\|_2$, and define $\mathcal A:H\to H^*$ by
\[
\langle \mathcal A v,\varphi\rangle
:=\int_\Omega\frac{v\varphi}{|x|^2}\,dx+(\nabla v,\nabla\varphi).
\]
Hardy's inequality gives
\[
\left|\int_\Omega\frac{v\varphi}{|x|^2}\,dx\right|
\le H_n\|v\|_H\|\varphi\|_H,
\]
while $\langle\mathcal Av,v\rangle\ge\|v\|_H^2$. Hence $\mathcal A$ is a bounded isomorphism from $H$ onto $H^*$ by the Lax--Milgram theorem.

For $t\ge0$ and $v\in H$, set
\[
\langle\mathcal F(t,v),\varphi\rangle
:=-(\nabla v,\nabla\varphi)
+k(t)\int_\Omega |v|^{p(x,t)-2}v\varphi\,dx.
\]
We verify the local Lipschitz property explicitly. Put
\[
2^*=\frac{2n}{n-2},\qquad q=(2^*)'=\frac{2n}{n+2},
\qquad r_*:=\frac{2^*}{p^+-1}.
\]
Since $p^+<2^*$, one has $r_*>q$. For $2\le r\le p^+$,
\[
\big||a|^{r-2}a-|b|^{r-2}b\big|
\le C_{p^+}(1+|a|+|b|)^{p^+-2}|a-b|.
\]
H\"older's inequality with exponents $2^*/(p^+-2)$ and $2^*$ therefore gives
\[
\big\||v|^{p(\cdot,t)-2}v-|w|^{p(\cdot,t)-2}w\big\|_{L^{r_*}(\Omega)}
\le C_R\|v-w\|_{L^{2^*}(\Omega)}
\]
whenever $\|v\|_H+\|w\|_H\le R$. Because $\Omega$ is bounded and $r_*>q$, the embedding $L^{r_*}(\Omega)\hookrightarrow L^q(\Omega)$ and the Sobolev embedding $H\hookrightarrow L^{2^*}(\Omega)$ imply
\[
\big\||v|^{p(\cdot,t)-2}v-|w|^{p(\cdot,t)-2}w\big\|_{L^q(\Omega)}
\le C_R\|v-w\|_H.
\]
Thus $v\mapsto\mathcal F(t,v)$ is locally Lipschitz from $H$ to $H^*$, uniformly for $t$ in compact intervals.

For fixed $v\in H$, continuity of $t\mapsto\mathcal F(t,v)$ follows from $p,k\in C^1$ and dominated convergence. Indeed,
\[
|v|^{p(x,t)-1}\le |v|+|v|^{p^+-1},
\]
and the right-hand side belongs to $L^q(\Omega)$ because $(p^+-1)q<2^*$. Consequently $(P)$ is equivalent to the nonautonomous ordinary differential equation
\[
u_t=\mathcal A^{-1}\mathcal F(t,u)
\]
in $H$. The standard local existence and continuation theorem for ordinary differential equations in Banach spaces \cite{Deimling} yields a unique maximal solution $u\in C^1([0,T^*);H)$, and a finite maximal time can occur only if $\|u(t)\|_H$ becomes unbounded.

Conversely, any weak solution in the sense above satisfies $\mathcal A u_t=\mathcal F(t,u)$ for a.e.\ $t$. Since $u\in C([0,T];H)$ and the right-hand side is continuous with values in $H^*$, $u_t$ has the continuous representative $\mathcal A^{-1}\mathcal F(t,u(t))$. Hence every weak solution coincides with the $C^1(H)$ solution supplied by the Banach-space ODE.

Finally, Hardy's inequality gives
\[
\frac12\|\nabla u(t)\|_2^2\le L(t)
\le\frac{1+H_n}{2}\|\nabla u(t)\|_2^2.
\]
Thus the continuation criterion is equivalent to the asserted divergence of $L$ along a sequence approaching $T^*$.
\end{proof}

\subsection{Energy estimates}

Let $p$, $k$ satisfy \eqref{eq:p-assumption} and \eqref{eq:k-assumption}, respectively, and let $T>0$. Recall the functionals $J$, $I$, and $E$ from \eqref{eq:J-def}--\eqref{eq:E-def}. The following identities are the basis of the blow-up arguments.

\begin{lemma}\label{lem:energy}
Let $u$ be a weak solution to $(P)$ on $[0,T)$. Define
\begin{equation}\label{eq:P-def}
\mathfrak P(t):=
\int_\Omega \frac{p_t(x,t)}{p(x,t)^2}
\bigl[p(x,t)\log |u(x,t)|-1\bigr]
|u(x,t)|^{p(x,t)}\,dx,
\end{equation}
where $|u|^{p}\log|u|$ is understood as $0$ at $u=0$. Then:
\begin{enumerate}
\item for every $t_0\in[0,T)$,
\begin{align}
J(u(t_0),t_0)
&+\int_0^{t_0}\Bigg[
\left\|\frac{u_t(s)}{|x|}\right\|_2^2
+\|\nabla u_t(s)\|_2^2 \notag\\
&\qquad
+k'(s)\int_\Omega\frac{|u(x,s)|^{p(x,s)}}{p(x,s)}\,dx
+k(s)\mathfrak P(s)
\Bigg]ds
=J(u_0,0);                                      \label{eq:energy-identity}
\end{align}
\item for a.e. $t\in(0,T)$,
\begin{equation}\label{eq:L-identity}
L'(t)
=\int_\Omega\frac{u(t)u_t(t)}{|x|^2}\,dx
+(\nabla u(t),\nabla u_t(t))
=-I(u(t),t).
\end{equation}
\end{enumerate}
\end{lemma}

\begin{proof}
Fix $T_0<T$ and define
\[
\Psi(t,v):=\int_\Omega\frac{|v(x)|^{p(x,t)}}{p(x,t)}\,dx,
\qquad 0\le t\le T_0,\quad v\in H:=W^{1,2}_0(\Omega).
\]
Since $p^+<2^*$, choose $\varepsilon>0$ such that $p^++\varepsilon<2^*$. The estimates
\[
s^r|\log s|\le C_\varepsilon(1+s^{r+\varepsilon}),
\qquad s\ge0,\quad p^-\le r\le p^+,
\]
and
\[
|s|^{r-1}\le |s|+|s|^{p^+-1}
\]
show, using the Sobolev embedding, that $\Psi$ is continuously differentiable along bounded subsets of $[0,T_0]\times H$. Its derivatives are
\[
D_v\Psi(t,v)[w]
=\int_\Omega |v|^{p(x,t)-2}vw\,dx
\]
and
\[
\partial_t\Psi(t,v)
=\int_\Omega\frac{p_t(x,t)}{p(x,t)^2}
\bigl[p(x,t)\log|v|-1\bigr]|v|^{p(x,t)}\,dx,
\]
with the logarithmic term interpreted as $0$ at $v=0$. To justify continuity of the second expression, let $(t_j,v_j)\to(t,v)$ in $[0,T_0]\times H$. Then $v_j\to v$ in $L^{p^++\varepsilon}(\Omega)$, while the preceding logarithmic growth estimate shows that the corresponding integrands are uniformly integrable. After passage to a subsequence they converge almost everywhere, and Vitali's theorem yields convergence in $L^1(\Omega)$. Thus $\partial_t\Psi$ is continuous. The same, simpler argument applies to $D_v\Psi$, and hence $\Psi\in C^1([0,T_0]\times H)$.

By Proposition \ref{prop:local}, a weak solution has a $C^1([0,T_0];H)$ representative. The ordinary Banach-space chain rule applied to $t\mapsto\Psi(t,u(t))$ therefore gives, for every $t\in[0,T_0]$,
\begin{align}
\frac{d}{dt}\int_\Omega\frac{|u|^{p(x,t)}}{p(x,t)}\,dx
={}&\int_\Omega |u|^{p(x,t)-2}u\,u_t\,dx \notag\\
&+\int_\Omega\frac{p_t(x,t)}{p(x,t)^2}
\bigl[p(x,t)\log|u|-1\bigr]|u|^{p(x,t)}\,dx.             \label{eq:variable-chain}
\end{align}
Since $T_0<T$ was arbitrary, this identity holds locally on $[0,T)$.

Using $u_t$ as a test function in \eqref{eq:weak-form} gives
\[
\left\|\frac{u_t}{|x|}\right\|_2^2
+\|\nabla u_t\|_2^2
+(\nabla u,\nabla u_t)
=k(t)(|u|^{p(x,t)-2}u,u_t).
\]
Combining this equality with \eqref{eq:variable-chain} yields, for a.e. $t\in(0,T)$,
\begin{align*}
\frac{d}{dt}J(u(t),t)
={}&-\left\|\frac{u_t(t)}{|x|}\right\|_2^2
-\|\nabla u_t(t)\|_2^2 \notag\\
&-k'(t)\int_\Omega\frac{|u(x,t)|^{p(x,t)}}{p(x,t)}\,dx
-k(t)\mathfrak P(t).
\end{align*}
Integration over $(0,t_0)$ proves \eqref{eq:energy-identity}.

Finally, differentiating \eqref{eq:L-def} and taking $\varphi=u(t)$ in \eqref{eq:weak-form} gives \eqref{eq:L-identity}.
\end{proof}

The next concavity argument is classic and is used extensively in the literature for a sufficient condition of blow-up time.

\begin{lemma} \label{lem:concavity}
	Let $\theta>0$ and let $\psi$ be positive and continuously differentiable on $[0,T)$, with $\psi'$ locally absolutely continuous, such that
	$\psi(0)>0$, $\psi'(0)>0$, and
	\[
	\psi(t)\psi''(t)-(1+\theta)(\psi'(t))^2\ge0,
	\qquad \text{for a.e. }0<t<T.
	\]
	Then
	\[
	t<\frac{\psi(0)}{\theta\psi'(0)}
	\qquad\text{for every }0<t<T.
	\]
	Consequently,
	\[
	T\le\frac{\psi(0)}{\theta\psi'(0)}.
	\]
\end{lemma}
\begin{proof}
	Set $y(t):=\psi(t)^{-\theta}$. Then $y'$ is locally absolutely continuous and, for a.e. $t$,
	\[
	y''(t)
	=-\theta\psi(t)^{-\theta-2}
	\left[\psi(t)\psi''(t)-(1+\theta)(\psi'(t))^2\right]\le0.
	\]
	Thus $y$ is concave and
	\[
	y(t)\le y(0)+t y'(0)
	=\psi(0)^{-\theta}
	\left(1-\theta\frac{\psi'(0)}{\psi(0)}t\right).
	\]
	Since $y(t)>0$, the factor in parentheses must be positive, which proves the claim.
\end{proof}

\section{Upper bounds for the blow-up time} \label{sec:upper}

In this section we work with the upper bounds for the blow-up time.
These are the contents of Theorems \ref{thm:negative} and \ref{thm:positive}.
Throughout this section, $L$ denotes the natural energy quantity defined in \eqref{eq:L-def}.

We start with the proof of Theorem \ref{thm:negative} which deals with negative modified initial energy.

\begin{proof}[Proof of Theorem \ref{thm:negative}]
	Let $T^*\in(0,\infty]$ be the maximal existence time of $u$.
	We show first that $T^*<\infty$ and then derive the asserted upper bound.

	Let $\mathfrak{P}$ be given by \eqref{eq:P-def}. For $s\ge0$, the function $s(1-\log s)$ (with value $0$ at $s=0$) is bounded above by $1$. Since
	$p(x,t)\log|u|=\log(|u|^{p(x,t)})$, $p_t\ge0$, and $k(t)\le k_\infty$, we have for a.e. $t\in(0,T^*)$,
	\begin{align}
	-k(t)\mathfrak{P}(t)
	&=-k(t)\int_\Omega \frac{p_t(x,t)}{p(x,t)^2}
	\big[p(x,t)\log|u(x,t)|-1\big]|u(x,t)|^{p(x,t)}\,dx \notag\\
	&\le k(t)\int_\Omega \frac{p_t(x,t)}{p(x,t)^2}\,dx \notag\\
	&\le k_\infty\int_\Omega \frac{p_t(x,t)}{p(x,t)^2}\,dx
	=-\frac{d}{dt}\int_\Omega \frac{k_\infty}{p(x,t)}\,dx.
	\label{eq:P-est}
	\end{align}

	Set
	\[
	K(t):=-E(u(t),t),\qquad 0\le t<T^*.
	\]
	By hypothesis, $K(0)=-E(u_0,0)>0$.
	Using Lemma \ref{lem:energy} together with \eqref{eq:P-est}, we obtain for a.e. $t\in(0,T^*)$,
	\begin{align}
	K'(t)
	={}&-\frac{d}{dt}\int_\Omega\frac{k_\infty}{p(x,t)}\,dx
	-\frac{d}{dt}J(u(t),t) \notag\\
	={}&\left\|\frac{u_t(t)}{|x|}\right\|_2^2
	+\|\nabla u_t(t)\|_2^2
	+k'(t)\int_\Omega\frac{|u(x,t)|^{p(x,t)}}{p(x,t)}\,dx \notag\\
	&\quad +k(t)\mathfrak{P}(t)
	-\frac{d}{dt}\int_\Omega\frac{k_\infty}{p(x,t)}\,dx \notag\\
	\ge{}&\left\|\frac{u_t(t)}{|x|}\right\|_2^2
	+\|\nabla u_t(t)\|_2^2\ge0.
	\label{eq:K-derivative}
	\end{align}
	Hence $K$ is nondecreasing and
	\[
	K(t)\ge K(0)>0,
	\qquad 0\le t<T^*.
	\]

	On the other hand, Lemma \ref{lem:energy} and $p(x,t)\ge p^-$ give
	\begin{align}
	L'(t)
	&=-I(u(t),t) \notag\\
	&=-\|\nabla u(t)\|_2^2
	+k(t)\int_\Omega |u(x,t)|^{p(x,t)}\,dx \notag\\
	&\ge \frac{p^- -2}{2}\|\nabla u(t)\|_2^2-p^-J(u(t),t) \notag\\
	&\ge -p^-J(u(t),t) \notag\\
	&=p^-\left[K(t)+\int_\Omega\frac{k_\infty}{p(x,t)}\,dx\right]
	\ge p^-K(t)>0
	\label{eq:L-prime-negative}
	\end{align}
	for a.e. $t\in(0,T^*)$.

	We next combine the two differential inequalities without discarding either component of the natural energy. By \eqref{eq:K-derivative},
	\begin{align*}
	L(t)K'(t)
	&\ge \frac12\left(
	\left\|\frac{u(t)}{|x|}\right\|_2^2+\|\nabla u(t)\|_2^2
	\right)
	\left(
	\left\|\frac{u_t(t)}{|x|}\right\|_2^2+\|\nabla u_t(t)\|_2^2
	\right)\\
	&\ge \frac12\left[
	\left(\frac{u(t)}{|x|^2},u_t(t)\right)
	+(\nabla u(t),\nabla u_t(t))
	\right]^2\\
	&=\frac12\big(L'(t)\big)^2
	\ge \frac{p^-}{2}K(t)L'(t),
	\end{align*}
	where the second inequality is the Cauchy--Schwarz inequality in
	$L^2(\Omega)\times L^2(\Omega;\mathbb{R}^n)$ and the last one follows from \eqref{eq:L-prime-negative}.
	Consequently,
	\[
	\left(K(t)L(t)^{-p^-/2}\right)'
	=L(t)^{-(p^-+2)/2}
	\left(L(t)K'(t)-\frac{p^-}{2}K(t)L'(t)\right)\ge0
	\]
	for a.e. $t\in(0,T^*)$.
	Thus
	\begin{equation}\label{eq:xi-lower}
	K(t)L(t)^{-p^-/2}\ge \xi_0
	:=K(0)L(0)^{-p^-/2}>0.
	\end{equation}
	Combining \eqref{eq:L-prime-negative} and \eqref{eq:xi-lower} yields
	\[
	L'(t)L(t)^{-p^-/2}\ge p^-\xi_0.
	\]
	Since $p^->2$,
	\[
	\frac{d}{dt}L(t)^{(2-p^-)/2}
	=\frac{2-p^-}{2}L(t)^{-p^-/2}L'(t)
	\le -\frac{p^-(p^- -2)}{2}\xi_0.
	\]
	Integrating over $(0,\tau)$, $0<\tau<T^*$, gives
	\begin{equation}\label{eq:negative-decay}
	0<L(\tau)^{(2-p^-)/2}
	\le L(0)^{(2-p^-)/2}
	-\frac{p^-(p^- -2)}{2}\xi_0\tau.
	\end{equation}
	The right-hand side of \eqref{eq:negative-decay} must remain positive for every $\tau<T^*$; hence
	\[
	T^*\le
	\frac{2L(0)^{(2-p^-)/2}}{p^-(p^- -2)\xi_0}
	=\frac{2L(0)}{p^-(p^- -2)K(0)}.
	\]
	Since $K(0)=-E(u_0,0)$ and
	\[
	2L(0)=\left\|\frac{u_0}{|x|}\right\|_2^2+\|\nabla u_0\|_2^2,
	\]
	we obtain exactly
	\[
	T^*\le
	\frac{\left\|u_0/|x|\right\|_2^2+\|\nabla u_0\|_2^2}
	{p^-(p^- -2)[-E(u_0,0)]}.
	\]
	In particular $T^*<\infty$. Proposition \ref{prop:local} now yields \eqref{eq:blowup}, and hence $u$ blows up at $T^*$.
\end{proof}

Next we prove Theorem \ref{thm:positive}, which treats nonnegative modified initial energy below the threshold in its statement.

\begin{proof}[Proof of Theorem \ref{thm:positive}]
	Let $T^*\in(0,\infty]$ be the maximal existence time of $u$, and write
	\[
	D(t):=\left\|\frac{u_t(t)}{|x|}\right\|_2^2+\|\nabla u_t(t)\|_2^2.
	\]
	The estimate \eqref{eq:K-derivative}, whose derivation does not use the sign of $E(u_0,0)$, yields
	\begin{equation}\label{eq:K-positive}
	K'(t)\ge D(t)\ge0,
	\qquad K(t):=-E(u(t),t).
	\end{equation}
	In particular,
	\begin{equation}\label{eq:K-integrated}
	K(t)\ge K(0)+\int_0^tD(s)\,ds
	=-E(u_0,0)+\int_0^tD(s)\,ds.
	\end{equation}

	By Hardy's inequality,
	\[
	\left\|\frac{u(t)}{|x|}\right\|_2^2\le H_n\|\nabla u(t)\|_2^2,
	\]
	and hence
	\begin{equation}\label{eq:gradient-controls-L}
	\|\nabla u(t)\|_2^2
	\ge \frac{2}{1+H_n}L(t).
	\end{equation}
	As in \eqref{eq:L-prime-negative}, but retaining the gradient contribution, we have
	\begin{align}
	L'(t)
	&\ge \frac{p^- -2}{2}\|\nabla u(t)\|_2^2-p^-J(u(t),t) \notag\\
	&=\frac{p^- -2}{2}\|\nabla u(t)\|_2^2
	+p^-K(t)+p^-\int_\Omega\frac{k_\infty}{p(x,t)}\,dx \notag\\
	&\ge \frac{p^- -2}{1+H_n}L(t)+p^-K(t). \label{eq:L-prime-positive}
	\end{align}
	Set
	\[
	a:=\frac{p^- -2}{1+H_n},
	\qquad
	C_1=\frac{p^-}{a}=\frac{p^-(1+H_n)}{p^- -2},
	\qquad
	M(t):=L(t)+C_1K(t).
	\]
	Then \eqref{eq:L-prime-positive} gives
	\[
	L'(t)\ge aM(t),
	\]
	while \eqref{eq:K-positive} implies
	\[
	M'(t)=L'(t)+C_1K'(t)\ge L'(t)\ge aM(t).
	\]
	Since
	\[
	M(0)=L(0)-C_1E(u_0,0)=M_0>0,
	\]
	Gronwall's inequality gives $M(t)\ge M_0e^{at}>0$. Therefore
	\begin{equation}\label{eq:L-increasing}
	L'(t)>0,
	\qquad L(t)\ge L(0),
	\qquad 0<t<T^*.
	\end{equation}

	We now apply a corrected concavity argument. Fix $0<\tau<T^*$ and choose
	\begin{equation}\label{eq:beta-choice}
	0<\beta<\frac{p^-M_0}{C_1(p^- -1)}.
	\end{equation}
	For a parameter $\sigma>0$ to be fixed below, define
	\[
	F(h):=\int_0^hL(s)\,ds+(\tau-h)L(0)+\frac{\beta}{2}(h+\sigma)^2,
	\qquad 0\le h\le\tau.
	\]
	Then $F(h)>0$ and, using $L'(s)=\langle u(s),u_t(s)\rangle_{\mathcal H}$ for the natural Hilbert product
	\[
	\langle v,w\rangle_{\mathcal H}
	:=\left(\frac{v}{|x|},\frac{w}{|x|}\right)_{L^2}
	+(\nabla v,\nabla w)_{L^2},
	\]
	we have
	\[
	F'(h)=\int_0^h\langle u(s),u_t(s)\rangle_{\mathcal H}\,ds
	+\beta(h+\sigma),
	\qquad
	F''(h)=L'(h)+\beta.
	\]
	Cauchy--Schwarz in $L^2(0,h;\mathcal H)\times\mathbb R$ gives
	\begin{align}
	(F'(h))^2
	&\le\left(2\int_0^hL(s)\,ds+\beta(h+\sigma)^2\right)
	\left(\int_0^hD(s)\,ds+\beta\right)\notag\\
	&\le 2F(h)\left(\int_0^hD(s)\,ds+\beta\right). \label{eq:F-CS}
	\end{align}

	Next, \eqref{eq:K-integrated} and $-J(u(h),h)=K(h)+\int_\Omega k_\infty/p(x,h)\,dx\ge K(h)$ yield
	\begin{align}
	L'(h)
	&\ge \frac{p^- -2}{2}\|\nabla u(h)\|_2^2+p^-K(h)\notag\\
	&\ge \frac{p^- -2}{2}\|\nabla u(h)\|_2^2
	-p^-E(u_0,0)+p^-\int_0^hD(s)\,ds. \label{eq:L-prime-concavity}
	\end{align}
	Put
	\[
	\theta:=\frac{p^- -2}{2}>0.
	\]
	Since $1+\theta=p^-/2$, inequalities \eqref{eq:F-CS} and \eqref{eq:L-prime-concavity} give
	\begin{align*}
	&F(h)F''(h)-(1+\theta)(F'(h))^2\\
	&\quad\ge F(h)\left[
	F''(h)-p^-\left(\int_0^hD(s)\,ds+\beta\right)
	\right]\\
	&\quad\ge F(h)\left[
	\frac{p^- -2}{2}\|\nabla u(h)\|_2^2
	-p^-E(u_0,0)-(p^- -1)\beta
	\right]\\
	&\quad\ge F(h)\left[
	\frac{p^- -2}{1+H_n}L(0)
	-p^-E(u_0,0)-(p^- -1)\beta
	\right]\\
	&\quad=F(h)\left[
	\frac{p^-}{C_1}M_0-(p^- -1)\beta
	\right]>0,
	\end{align*}
	where we used \eqref{eq:gradient-controls-L}, \eqref{eq:L-increasing}, and \eqref{eq:beta-choice}.

	Lemma \ref{lem:concavity} therefore yields
	\[
	\tau\le\frac{F(0)}{\theta F'(0)}
	=\frac{2\tau L(0)+\beta\sigma^2}
	{(p^- -2)\beta\sigma}.
	\]
	Choose
	\begin{equation}\label{eq:sigma-choice}
	\sigma>\frac{2L(0)}{(p^- -2)\beta}.
	\end{equation}
	Then
	\[
	\tau\le
	\frac{\beta\sigma^2}
	{(p^- -2)\beta\sigma-2L(0)}.
	\]
	The right-hand side is minimized at
	\[
	\sigma=\frac{4L(0)}{(p^- -2)\beta},
	\]
	which is admissible in \eqref{eq:sigma-choice}, and hence
	\begin{equation}\label{eq:tau-beta}
	\tau\le\frac{8L(0)}{(p^- -2)^2\beta}.
	\end{equation}
	Because \eqref{eq:tau-beta} holds for every $\beta$ satisfying \eqref{eq:beta-choice}, letting
	\[
	\beta\uparrow\frac{p^-M_0}{C_1(p^- -1)}
	\]
	gives
	\[
	\tau\le
	\frac{8C_1(p^- -1)L(0)}
	{p^-(p^- -2)^2M_0}.
	\]
	Finally let $\tau\uparrow T^*$. We obtain the asserted finite upper bound for $T^*$.
	Proposition \ref{prop:local} then gives \eqref{eq:blowup}, so the solution blows up at $T^*$.
\end{proof}

\section{A lower bound for the blow-up time} \label{sec:lower}

In this section we provide a lower bound for the blow-up time, which is Theorem \ref{thm:lower}.
We continue to use the natural energy quantity $L$ defined in \eqref{eq:L-def}.

\begin{proof}[Proof of Theorem \ref{thm:lower}]
Let $T^*$ be the blow-up time. For every $h\in[0,T^*)$, the definitions of $p^-$ and $p^+$ give the pointwise estimate
\begin{equation}\label{eq:variable-power}
|u(x,h)|^{p(x,h)}\le |u(x,h)|^{p^-}+|u(x,h)|^{p^+}.
\end{equation}
Indeed, the first term on the right dominates on $\{|u|\le1\}$, while the second dominates on $\{|u|>1\}$.

Since the standing assumption \eqref{eq:p-assumption} implies $2<p^-\le p^+<2^*$, Lemma \ref{lem:sobolev} applies with $q=p^-$ and $q=p^+$. Thus, using \eqref{eq:I-def}, \eqref{eq:variable-power}, and $k(h)\le k_\infty$, we obtain
\begin{align}
L'(h)
&=-I(u(h),h)\notag\\
&=k(h)\int_\Omega |u(x,h)|^{p(x,h)}\,dx
-\|\nabla u(h)\|_2^2\notag\\
&\le k_\infty S_{p^-}^{p^-}\|\nabla u(h)\|_2^{p^-}
+k_\infty S_{p^+}^{p^+}\|\nabla u(h)\|_2^{p^+}\notag\\
&\le C^*\left(L(h)^{p^-/2}+L(h)^{p^+/2}\right),
\label{eq:lower-ode}
\end{align}
where one may take
\[
C^*:=k_\infty\max\left\{
2^{p^-/2}S_{p^-}^{p^-},
2^{p^+/2}S_{p^+}^{p^+}
\right\}.
\]
In the last step we used $\|\nabla u(h)\|_2^2\le2L(h)$.
With $\gamma^\pm=p^\pm/2$, estimate \eqref{eq:lower-ode} becomes
\[
L'(h)\le C^*\bigl(L(h)^{\gamma^-}+L(h)^{\gamma^+}\bigr).
\]
Notice that $\gamma^-\le\gamma^+$ and $\gamma^->1$.

Because $u$ is a nontrivial solution which blows up at $T^*$, uniqueness in Proposition \ref{prop:local} implies $L(h)>0$ for every $h\in[0,T^*)$. Fix $t_0\in[0,T^*)$ and define
\[
\Phi(r):=\int_{L(t_0)}^r
\frac{ds}{s^{\gamma^-}+s^{\gamma^+}},
\qquad r>0.
\]
The chain rule and the preceding differential inequality give, for a.e. $h\in(t_0,T^*)$,
\[
\frac{d}{dh}\Phi(L(h))
=\frac{L'(h)}{L(h)^{\gamma^-}+L(h)^{\gamma^+}}
\le C^*.
\]
Consequently,
\[
\int_{L(t_0)}^{L(t)}
\frac{ds}{s^{\gamma^-}+s^{\gamma^+}}
\le C^*(t-t_0),
\qquad t_0<t<T^*.
\]
Because $u$ blows up at $T^*$, there exists a sequence $t_j\uparrow T^*$ such that $L(t_j)\to\infty$. Since $\gamma^->1$, the improper integral converges. Letting $j\to\infty$ yields
\[
T^*\ge t_0+\frac1{C^*}
\int_{L(t_0)}^\infty
\frac{ds}{s^{\gamma^-}+s^{\gamma^+}}.
\]

If $L(t_0)\ge1$, then $s^{\gamma^-}\le s^{\gamma^+}$ for every $s\ge L(t_0)$, and hence
\[
\int_{L(t_0)}^\infty
\frac{ds}{s^{\gamma^-}+s^{\gamma^+}}
\ge\frac12\int_{L(t_0)}^\infty s^{-\gamma^+}\,ds
=\frac{L(t_0)^{1-\gamma^+}}{2(\gamma^+-1)}.
\]
This proves both asserted estimates.
\end{proof}

\section*{Statements and Declarations}

\textbf{Competing interests.} The authors declare that they have no competing interests.

\textbf{Data availability.} No datasets were generated or analyzed in this theoretical study.

\bibliographystyle{amsplain}

\end{document}